\documentclass[reqno,fleqn,11pt]{amsart}

\usepackage{amsmath}
\usepackage{amssymb}
\usepackage{amsfonts}
\usepackage{amsthm}

\newcommand{\C}{\mathbb{C}}
\newcommand{\f}{\rightarrow}
\newcommand{\deb}{\bar\partial}
\newcommand{\de}{\partial}
\newcommand{\K}{K\"{a}hler}

\newcommand{\us}{\underline s}

\newcommand{\GL}{\operatorname{GL}}
\newcommand{\Hom}{\operatorname{Hom}}

\newcommand{\rank}{\operatorname{rank}}
\newcommand{\diag}{\operatorname{diag}}
\newcommand{\Ric}{\operatorname{Ric}}

\newtheorem{prop}{Proposition}
\newtheorem{thm}[prop]{Theorem}

\newtheorem{cor}[prop]{Corollary}

\newcommand{\fdim}{\hspace*{\fill}$\Box$}

\begin{document}

\title{Balanced metrics on homogeneous vector bundles}
\author[R. Mossa]{Roberto Mossa}
\address{Dipartimento di Matematica e Informatica, Universit\`{a} di Cagliari,
Via Ospedale 72, 09124 Cagliari, Italy}
\email{roberto.mossa@gmail.com}
\thanks{Research partially supported by GNSAGA (INdAM) and MIUR of Italy}
\date{January 15, 2011}

\begin{abstract}
Let $E\rightarrow M$ be a holomorphic vector bundle over a compact \K\ manifold $(M, \omega)$
and let $E=E_1\oplus\cdots \oplus E_m\rightarrow M$ be its decomposition into  irreducible factors.
Suppose that each $E_j$ admits a $\omega$-balanced metric in Donaldson-Wang terminology.
In this paper we prove that $E$ admits a unique $\omega$-balanced metric if and only if $\frac{r_j}{N_j}=\frac{r_k}{N_k}$ for all  $j, k=1, \dots, m$,
where $r_j$ denotes the rank of $E_j$ and $N_j=\dim H^0(M, E_j)$.
We apply our result to the case of homogeneous vector bundles over a  rational homogeneous variety $(M, \omega)$ and  we  show the existence and rigidity of balanced \K\ embedding from $(M, \omega)$  into Grassmannians.

\end{abstract}

\maketitle

\section{Introduction}\label{sectionwang}

Let $E\rightarrow M$ be a very ample  holomorphic vector bundle of $\rank$ $r$ with $\dim (H^0(M,E))=N$. Choose a basis $\underline s=(s_1, \dots , s_N)$ of $H^0(M, E)$,  the space of  global holomorphic sections of $E$, we denotes with $i_{\underline s}:M\rightarrow G(r, N)$ the
Kodaira map associated to  the basis $\underline s$ (see, e.g. \cite{loimossa}).

Consider the flat metric  $h_0$ on the tautological bundle $\mathcal T \rightarrow G(r, N)$,
i.e.  $h_0(v,w)=w^*v$,  and the
 dual metric  $h_{Gr}=h_0^*$ on the quotient bundle $\mathcal Q=\mathcal T^*$. 
Let  $\omega_{Gr}=P^*\omega_{FS}$ be the \K\ form on $G(r, N)$  pull-back of the Fubini--Study form  $\omega_{FS}=$ \ $\frac{i}{2} \de \deb \log (|z_0|^2+\dots+|z_{\widetilde N-1}|^2)$ via  the Pl\"ucker embedding $P:G(r,N)\f \C P ^{\widetilde N -1}$, with $\widetilde N={\binom{N}{r}}$.

 We can endow
 $E=i_{\underline s}^*\mathcal Q$ with the hermitian   metric
\begin{equation}\label{defhs}
h_{\underline s}=i_{\underline s}^*h_{Gr}
\end{equation}
 and define a   $L^2$-product on $H^0(M,E)$
by the formula:
\begin{equation}\label{deflhs}
\langle \cdot, \cdot \rangle_{{h_{\underline s}}} =\frac{1}{V(M)} \int_M h_{\underline s}(\cdot, \cdot) \frac{\omega^n}{n!},
\end{equation}
where $\omega^n=\omega\wedge\cdots\wedge\omega$ and $V(M) = \int_M \frac{\omega^n}{n!}$.

An hermitian  metric $h$ on $E$
 is called $\omega$-{\em balanced} if there exists  a basis $\underline s$
 of $H^0(M, E)$
such that $h=h_{\underline s}=i_{\underline s}^*h_{Gr}$ and
\begin{equation}\label{balancedcondition}
\langle s_j, s_k\rangle_{{h_{\underline s}}}=\frac{r}{N} \delta_{jk}, \quad j, k=1, \dots , N=\dim H^{0}(M, E).
\end{equation}

\vskip 0.3cm
The concept of balanced metrics on complex vector bundles was introduced by X. Wang \cite{wang1} (see also \cite{wang2}) following S. Donaldson’s ideas introduced in \cite{donaldson} in order to characterize which manifolds admit a constant scalar curvature \K\ metric. While existence of $\omega$-balanced metrics is still a very difficult and obscure problem, uniqueness had been proved by Loi--Mossa (cfr. Theorem \ref{unicity}). Study of balanced metrics is a very fruitful area of research both from mathematical and physical point of view (see, e.g.,  \cite{cgr3}, \cite{cgr4}, \cite{culoi},  \cite{mebal}, \cite{baldisk}, \cite{regscal}  and  \cite{albergbal}).

Regarding the existence and uniqueness  of $\omega$-balanced basis we have the following fundamental results (see next section for the definition
of the Gieseker point $T_E$ of $E$).

\begin{thm}\label{wbal}(X. Wang, \cite{wang2}) The Gieseker point $T_E$ is stable (in the GIT terminology) if and only if $E$ admits a $\omega$-balanced metric.
\end{thm}

\begin{thm}\label{unicity}(A. Loi--R. Mossa, \cite{loimossa})
Let $E$ be a holomorphic vector bundle over a compact \K\ manifold
$(M, \omega)$. If  $E$ admits a  $\omega$-balanced metric  then the metric is unique.
\end{thm}

The following theorem is the main result of this paper  (the \lq\lq only if''  part is already contained  in Theorem 2 of \cite{loimossa} and we include it here for completeness).

\begin{thm}\label{splitbal}
Let $(M, \omega )$ a \K\ manifold and let $E_j \f M$ be a very ample vector bundles over $M$ with $\rank E_j=r_j$ and $\dim H^0(M,E_j)=N_j> 0$, for $j=1, \dots, m$. Suppose that each $E_j$ admits a $\omega$-balanced metric.
Then the vector bundle $E=\oplus _{j=1}^m E_j\rightarrow M$ admits a unique $\omega$-balanced metric if and only if $\frac{r_j}{N_j}=\frac{r_k}{N_k}$ for all  $j, k=1, \dots, m$.
\end{thm}

When $(M, \omega )$ is a rational homogeneous variety and $E \f M$ is an irreducible homogeneous vector bundle over $M$, L. Biliotti and A. Ghigi  \cite{ghigi} proved that  $E$  admits a unique $\omega$-balanced metric. Hence we immediately  get the following result.

\begin{cor}\label{rathom}
Let $(M, \omega )$ be a rational homogeneous variety and let $E_j \f M$ be  irreducible homogeneous vector bundles over $M$ with $\rank E_j=r_j$ and $\dim H^0(M,E_j)=N_j> 0$, for $j=1, \dots, m$.
Then the homogeneous vector bundle $E=\oplus _{j=1}^m E_j\rightarrow M$ admits a unique homogeneous $\omega$-balanced metric if and only if $\frac{r_j}{N_j}=\frac{r_k}{N_k}$ for all  $j, k=1, \dots, m$.
\end{cor}

The paper contained other two sections. In the next one we prove Theorem \ref{splitbal}. In the last section we apply our result to prove the existence and rigidity of balanced maps of rational homogeneous varieties into grassmannians.

\vskip 0.3cm

\noindent
{\bf Acknowledgements:} I wish to thank Prof. Alessandro Ghigi and Prof. Andrea Loi for various interesting and stimulating discussions.

\section{Proof of Theorem \ref{splitbal}}
The  Gieseker point $T_E$ of a vector bundle $E$ of rank $r$ is the map
\[
T_E:\bigwedge^r H^0(M,E) \f H^0(M, \det E)
\]
which sends  $s_1 \wedge \dots \wedge s_r \in \bigwedge^r H^0(M,E)$ to the holomorphic section  of $\det E$ defined by
\[
T_E(s_1 \wedge \dots \wedge s_r) :x\mapsto s_1(x) \wedge \dots \wedge s_r(x).
\]

The group $\GL(H^0(M, E))$ acts on $H^0(M,E)$, therefore we get an action also on $\bigwedge^r H^0(M,E)$ and on $\Hom(\bigwedge^r H^0(M,E), H^0(M,\det E))$. The actions
are given by
\[
V \cdot (s_1 \wedge \dots \wedge s_r)=Vs_1 \wedge \dots \wedge Vs_r
\]
and
\[
(V\cdot T)(s_1 \wedge \dots \wedge s_r)=T(V \cdot (s_1 \wedge \dots \wedge s_r)).
\]
where $T \in \Hom\left(\bigwedge^r H^0(M,E),H^0(M,\det E)\right)$ and $V \in \GL\left(H^0(M, E)\right)$.

Recall also that  if $G$ is a reductive group acting linearly on a vector space V then 
an element $v$ of $V$ is called \emph{stable} if $Gv$ is closed in $V$ and the stabilizer of $v$ inside $G$ is finite.

\vskip 0.3cm

\noindent
{\bf Proof of Theorem \ref{splitbal}} Without loss of generality we can assume $m=2$.
Hence assume that $E$ is a direct sum of two holomorphic vector bundles $E_1, E_2 \f M$ with $\rank E_j=r_j$ and $\dim H^0(M, E_j)=N_j>0$, $j=1,2$. 

Suppose first   that $\frac{N_1}{r_1}= \frac{N_2}{r_2}$. Let ${\underline s}^j=\{s^j_1, \dots, s^j_{N_j}\}$ be the basis of $H^0(M,E_j)$ for $j=1, 2$.
Then, the assumption $\frac{r_1}{N_1}=\frac{r_1}{N_2}$,
readily implies that  the basis
\begin{equation}\label{beq}
\underline s=((s^1_1,0), \dots, (s_{N_1}^1,0),(0,s_1^2),\dots, (0,s_{N_2}^2))
\end{equation}
is a homogeneous $\omega$-balanced basis for $E_1 \oplus E_2$. Therefore $h_{\underline s}=i_{\underline s}^*h_{Gr}$ is the desired
homogeneous balanced metric on
$E_1 \oplus E_2$ which is unique   by Theorem \ref{unicity}.

\vskip 0.cm

Viceversa  if $\frac{N_1}{r_1}\neq \frac{N_2}{r_2}$ we claim that  the Gieseker point of $E$ is not stable
and by  Wang's Theorem \ref{wbal} the bundle $E$ can not admit a $\omega$-balanced metric.
In order to prove this consider the basis $\us = \{s_1,\dots,s_{N_1 + N_2}\}$ of $H^0(M,E)$ such that $\{s_1, \dots, s_{N_1}\}$ is a basis of $H^0(M,E_1\oplus \{0\})$ and $\{s_{N_1+1},\dots,s_{N_1+N_2}\}$ is a basis of $H^0(M,\{0\} \oplus E_2)$. Suppose that $\frac{N_1}{r_1} > \frac{N_2}{r_2}$.
Consider the $1$-parameter subgroup of $SL(N_1+N_2)$
\[
t \mapsto g(t)= \diag(\underbrace{t^{-N_2},\dots, t^{-N_2}}_{\text{$N_1$ times}},t^{N_1},\dots,t^{N_1}),
\]
where the action on the elements of the basis $\us$ is
\[
g(t)s_j=
\left\{ \begin{array}{ll}
  t^{-N_2} s_j & \quad \text{if $ j \leq N_1$} \\
  t^{N_1}s_j & \quad \text{otherwise.}\\
\end{array} \right.
\]
Observe that the section $x \mapsto s_{j_1}(x)\wedge \dots \wedge s_{j_r}(x)$ (with $j_1<j_2<\dots<j_r$), were $r=r_1+r_2$, is different from zero only when $j_{r_1}\leq N_1 < j_{r_1+1}$. So the action of $g(t)$ on the Gieseker point is given by:
\[
g(t)T_E(s_{j_1}\wedge \dots \wedge s_{j_r})=
\]
\[
=\left\{
\begin{array}{l l}
  t^{r_2N_1-r_1N_2}T_E(s_{j_1}\wedge \dots \wedge s_{j_r}) & \quad \text{if $j_{r_1}\leq N_1<j_{r_1+1}$} \\
  0 & \quad \text{otherwise}\\
\end{array} \right.
\]
Since $\frac{N_1}{r_1} > \frac{N_2}{r_2}$ this yields
\[
\lim_{t \f 0} g(t) T_E \equiv 0.
\]
and hence $T_E$ is not stable with respect to the reductive group $\GL (H^0(M, E)=\GL (N_1+N_2)$.
\fdim

\section{Balanced maps of rational homogeneous varieties into grassmannians}
Given a \K\ manifold $(M,\omega)$ we call $f:M \f G(r,N)$ a balanced \K\ embedding  if it is a \K\ embedding with respect to the \K\ form $\omega_{Gr}$ on $G(r,N)$ and it is of the form $f=i_{\us}$, where $\us$ is a $\omega$-balanced basis for $H^0(M,E)$ (i.e. the condition \eqref{balancedcondition} is satisfied), with $E=f^*\mathcal Q$.  

\begin{thm}\label{rigidity}
Let $(M,\omega)$ be a rational homogeneous variety and let $E \f M$ be an homogeneous vector bundle over $M$ with $\omega\in c_1(E)$. Suppose that $E$ is direct sum of irreducible homogeneous vector bundles $E_1,\dots,E_m \f M$ with $\rank E_j=r_j$ and $\dim H^0(M,E_j)=N_j> 0$, $j=1, \dots, m$. Then  there exists a unique (up to a unitary transformations of $G(r, N))$
balanced \K\ embedding $f:M\rightarrow G(r, N)$ such that
$f^*Q=E$ if and only if $\frac{r_j}{N_j}=\frac{r_k}{N_k}$ for all  $j, k=1, \dots, m$.
\end{thm}

\proof
Theorem \ref{rathom} tell us that there exists a balanced embedding if and only if is satisfied the condition $\frac{r_1}{N_1}=\dots=\frac{r_m}{N_m}$. Suppose that $f:M \f G(r,N)$ is the balanced \K\ embedding. Let $g$ an element of $G$ the group of isometries of $M$, then the map $f \circ g: M \f G(r,N)$ is again a balanced \K\ embedding, so by Theorem \ref{unicity} we get $(f\circ g)^* h_{Gr}=f^* h_{Gr}$. Therefore the \K\ form $i^*\omega_{Gr}$ is homogeneous, indeed
\[g^*(f^* \omega_{Gr}) = (f\circ g)^* \omega_{Gr}= \Ric( (f\circ g)^* h_{Gr})=\Ric( f^* h_{Gr})=f^*\omega_{Gr}.\]
Since over a rational homogeneous variety every cohomology class contains exactly one invariant form, then $i^*\omega_{Gr}=\omega$. 
\endproof

As a corollary of the previous result we get: 

\begin{cor}
Let $(M,\omega)$ be an hermitian symmetric space of compact type (HSSCT) of dimension $r$ and with  $\dim H^0(M,TM)=N$. Assume that $(M,g)$ is the product of irreducible  hermitian symmetric spaces of compact type $(M_1,\omega_1), \dots , (M_m,\omega_m)$. Then $(M,\omega)$ admit a balanced \K\ embedding in $G(r,N)$ if and only if   $\frac{r_j}{N_j}=\frac{r_k}{N_k}$ for all  $j, k=1, \dots, m$, where $\dim (M_j)=r_j$ and $\dim H^0(M,TM_j)=N_j> 0$.
\end{cor}

\end{document}